\documentclass[12pt]{iopart}

\usepackage{enumerate}
\usepackage[normalem]{ulem}
\usepackage{graphicx}
 \usepackage{subfigure}
 
\usepackage{lipsum}
\usepackage{graphicx}
\usepackage{epstopdf}
\usepackage{enumerate}
\usepackage{algorithm}
\newcommand{\bH}{{\bf H}}

\newcommand{\bD}{{\bf D}}

\newcommand{\bm}{{\bf m}}

\newcommand{\bW}{{\bf W}}
\newcommand{\bw}{{\bf w}}

\newcommand{\bui}{{\bf u_i}}
\newcommand{\bQ}{{\bf Q}}

\newcommand{\bqi}{{\bf q_i}}
\newcommand{\bPI}{{\bf P_i}}
\newcommand{\bP}{{\bf P}}

\newcommand{\bF}{{\bf F}}
\newcommand{\bB}{{\bf B}}
\newcommand{\bG}{{\bf G}}
\newcommand{\bM}{{\bf M}}
\newcommand{\bg}{{\bf g}}
\newcommand{\bDf}{{\bf D_f}}
\newcommand{\bDs}{{\bf D_s}}

\newcommand{\bL}{{\bf L}}
\newcommand{\bR}{{\bf R}}

\usepackage{algorithm}
\usepackage{algorithmicx}
\usepackage[noend]{algpseudocode}
\newcommand*\Let[2]{\State #1 $\gets$ #2}
\algrenewcommand\algorithmicrequire{\mboxbf{Input:}}
\algrenewcommand\algorithmicensure{\mboxbf{Output:}}

\usepackage{iopams}  
\begin{document}

\title[Inversion and Interpolation]{Simultaneous shot inversion for nonuniform geometries using fast data interpolation.}

\author{Michelle Liu\footnote{Department of Mathematics, University of British Columbia, Vancouver, BC, Canada}, 
Rajiv Kumar\footnote{School of Earth and Atmospheric Sciences, Georgia Institute of Technology, USA} 
Eldad Haber\footnote{Department of Earth and Ocean Sciences, University of British Columbia, Vancouver, BC, Canada}, 
and Aleksandr Aravkin\footnote{Department of Applied Mathematics, University of Washington, Seattle, WA, USA}.}

%\address{IOP Publishing, Temple Circus, Temple Way, Bristol BS1 6HG, UK}
%\ead{submissions@iop.org}
%\vspace{10pt}
%\begin{indented}
%\item[]August 2017
%\end{indented}

\begin{abstract}
Stochastic optimization is key to efficient inversion in PDE-constrained optimization. 
Using `simultaneous shots', or random superposition of source terms, works very well 
in simple acquisition geometries where all sources see all receivers, but this rarely occurs in practice. 

We develop an approach that interpolates data to an ideal acquisition geometry while 
solving the inverse problem using simultaneous shots. The approach is formulated as a joint 
inverse problem, combining ideas from low-rank interpolation with  full-waveform 
inversion. Results using synthetic experiments illustrate the flexibility and efficiency of the approach.
\end{abstract}

\vspace{2pc}
\noindent{\it Keywords}: Optimization, low-rank interpolation, full-waveform inversion. 
%\maketitle

\section{Introduction}
 Large scale inverse problems with partial differential equation (PDE) constraints play a key role 
  in many applications, including medical, electromagnetic, seismic imaging, as well as DC resistivity and hydrogeology~\cite{yaoguo1994inversion,steklova2017computational,pratt1999seismic,cheney1999electrical,haber2004inversion}.
We focus on seismic data, which are used by both global seismologists and oil and gas industries to get subsurface information of the Earth. 
In the exploratory setting of marine acquisition, seismic data are obtained by a ship towing a compressed air gun (i.e. a seismic source) and a stream of receivers. The air gun produces acoustic wave that propagate deep into the ocean floor, where a part of the wave is then get reflected to the surface where it gets recorded by the receivers. Given the observed data recorded at receiver, we solve for subsurface properties of earth such as velocity, density, and conductivity~\cite{virieux2009overview}.

Waveform inversion for the medium parameters can be formalized as an inverse problem with PDE constraints: 
\begin{eqnarray}
\label{eq:expanded}
&\min_{{\bf{u_1, . . ., u_N, m}}} \qquad \frac{1}{2N} \sum_{i=1}^N \| {\bf{P_i u_i - d_i}}\|_2^2\nonumber\\
& \mbox{subject to} \qquad {\bf{c_i}}({\bf{u_i, m}}):= {{\bH}(\bf{m})u_i - q_i} = 0, \;\>\>\>\>\> i = 1, \ldots, N,
\end{eqnarray}

\noindent where the constraints ${\bf{c_i(u_i}}, {\bf{m}})$ are the discretized linear PDEs, $N$ is the number of sources or experiments in a given survey, ${\bf{q_i}}\in\mathbb{R}^{l_q}$ represents the $i^{th}$ source that emits the field ${\bf{u_i}}\in\mathbb{R}^{l_u}$ , ${\bf{P_i}} $ is the matrix that maps the discretized field ${\bf{u_i}}$ to the location of where the data ${\bf{d_i}} \in \mathbb{R}^l$  was collected, and the matrix ${\bf{H}}$ is a discretization of the PDE with appropriate boundary conditions. We assume it is possible to compute the field ${\bf{u_i}}$ given {\bf{m}}:
\begin{equation}
\label{eq:ui}
\bui = \bH(\bm)^{-1}\bqi,
\end{equation}
and then (\ref{eq:expanded}) can be written in its {\it reduced} form:
 \begin{equation}
 \label{eq:reduced}
 \min_{\bm}\sum_{i=1}^N \|{\bPI\bH^{-1}}({\bf{m})q_i - d_i}\|_2^2 = \|\bP \bH^{-1}\bQ - \bD\|_F^2.
\end{equation}
  where ${\bf{F_i({m})}}:= {\bf{P_iH^{-1}}}({\bf{m})q_i}$ is the forward problem that predicts the data set for $i^{th}$ source.  
Prior information can also be included in this formulation using constraints or regularization for the model $\bm$, but
we do not focus on this in the paper.

 In a realistic setting, we have to solve large number of PDEs to evaluate $\bF$, with $N$ easily of the orders of $10^4$. 
Since each evaluation requires solving the linear system of equations~(\ref{eq:reduced}), 
these solves are the main computational bottleneck.  
Problem~(\ref{eq:reduced}) is written naturally as a large sum, so stochastic techniques can readily apply.
In particular, one can sample shots and use these smaller samples to generate updates for $\bm$~\cite{li2012fast}.
However, problem~(\ref{eq:reduced}) has additional structure that allows a specialized randomized approach. 
 
The simultaneous source method~\cite{van2013fast,roosta2015randomized,HCH} uses the following simple fact: 
\begin{equation}
\label{eq:ss}
\mathbb{E}_{\bw} \|\bP \bH^{-1}\bQ \bw - \bD\bw\|_F^2  =  \|\bP \bH^{-1}\bQ - \bD\|_F^2,
\end{equation}
for any $\bw$ satisfying $\mathbb{E}_{\bw} \bw\bw^T = I$. This makes it possible to create 
so called {\it simultaneous shots} $\widetilde \bqi  = \bQ\bw$ during the optimization process,
and obtain updates in $\bm$ using these shots rather than the entire dataset $\bQ$.  
Unfortunately, the identity~(\ref{eq:ss}) pre-supposes a common acquisition domain for 
all sources; it assumes that all sources see all receivers. This assumption is routinely violated
for many types of data acquisition, and in particular by the marine acquisition scenario.

\noindent
{\bf Contribution:}
The current paper makes it possible to use simultaneous shots in complex geometries, 
by treating uncollected data as if {\it missing} from a full all-see-all acquisition scenario, 
and using interpolation techniques to fill it in.  Given a dataset collected from a particular 
acquisition, we can use low rank regularization to interpolate the unobserved data, 
and then proceed with simultaneous sources. However, we also go further: 
we propose and solve a {\it joint inversion and interpolation} problem, which iteratively 
refines the interpolated data as the model estimates improve. 

\noindent
{\bf Roadmap:} The paper proceeds as follows. In Section~\ref{sec:background}, we review the simultaneous shots method and its relationship to trace estimation and stochastic optimization. 
We also formalize why the approach fails in complex acquisition geometries.  
In Section~\ref{sec:lowrank} we discuss low-rank interpolation, and show 
the efficacy and limitations of a two stage approach; (1) interpolation followed by 
(2) inversion by simultaneous shots. In Section~\ref{sec:unified}, we develop
a unified formulation that solves a single problem to accomplish the tasks simultaneously,
and show that this approach significantly improves the results. We also develop an algorithm for this 
unified optimization problem.  
We illustrate the advantages of the approach by inverting the Marmousi model \cite{brougois1990marmousi} 
in a difficult acquisition scenario. The unified approach gives much 
better results than the simple stage-wise workflow. 

\section{Background} \label{sec:background}
 
Suppose that we have an all-see-all geometry, so all sources are recorded by all the receivers,  
$\bP = \bPI, \forall i$. 
Define 
 \begin{equation}
\label{1.3}
\phi(\bm) := \frac{1}{2}\| \bP \bH^{-1}(\bm)\bQ - \bD \|^2_F
\end{equation} 
\noindent where ${\bf{Q = [q_1, . . . q_i, . . . ,q_N]}}$,  ${\bf{D = [d_1, . . . , d_i, . . . ,d_N]}}$ and ${\bf{q_i}}, {\bf{d_i}} \in \mathbb{R}^{l}$.  Then the identity~(\ref{eq:ss}) is easily derived for  $\bB = \bP \bH^{-1}(\bm)\bQ - \bD$:
 \begin{equation}
\mathbb{E}\|\bB \bw\|^2 = \mathbb{E} ( \bw^T \bB^T \bB \bw)  = \mathrm{tr}(\bB^T\bB) \mathbb{E} (\bw \bw^T) = \|\bB\|_F^2
\end{equation}
where the last equality assumes $\mathbb{E} (\bw \bw^T) = I$, which is true e.g. for a standard Gaussian vector. 
Many other distributions, e.g. Rademacher, also satisfy this requirement~\cite{hutchinson1990stochastic,haber2012effective}.

If we now consider the randomized function 
\[
\phi_\bw(\bm) =  \frac{1}{2}\| \bP \bH^{-1}(\bm)\bQ\bw - \bD\bw \|^2_2,
\]
we have $\mathbb{E}\phi_\bw(\bm)  = \phi(\bm)$, and evaluating $\phi_\bw(\bm)$ requires solving a single PDE. 
We can think of the original problem as a stochastic optimization problem 
\begin{equation}
\label{eq:stochopt}
\min_{\bm} \mathbb{E}( \| ({\bf{PH(m)^{-1} Q - D)w}}\|^2).
\end{equation}
This leads to two views: stochastic average approximation, and stochastic optimization.  

{\bf Stochastic Average Approximation (SAA)\cite{shapiro2009lectures} \cite{shapiro2003monte}.}
We think of approximating the expectation  by a small finite sample: 
 \begin{equation}
\label{eq:wtrick}
\mathbb{E} ( {\bf{w^T B^T B w}}) \approx  \hat{\phi}_K(\bf{m})
:= \frac{1}{K} \sum_{j=1}^K \| \bf{B}\bf{w}_j\|_2^2 
= \frac{1}{K} \sum_{j=1}^K \|{\bf{ PH(m)^{-1}(Qw_j) - Dw_j}}\|^2 
\end{equation}
\noindent where ${\bf{w_j}} \in\>\mathbb{R}^N$ is any random vector that satisfies $\mathbb{E}(\bw\bw^T = I)$. 
The solution to~(\ref{eq:wtrick}) is then approximated by the solution to the following problem:
  \begin{equation}
\label{eq:saa}
\min_\bm\> \frac{1}{2K} \sum_{j=1}^{K} ({\bf{ PH(m)^{-1}(Qw_j) - Dw_j}} )^2.
\end{equation}
We work with K PDEs, reducing the computational cost substantially if  $K<<N$. 

{\bf Stochastic optimization}. At any iteration, we can get an unbiased estimate of $\phi(\bm)$ and 
its gradient. Assuming for simplicity that we do not have any regularizer, 
we get a family of algorithms of the form 
\[
\bm^+ = \bm - \bG(\bm)^{-1} \bg(\bm), 
\]
where 
\[
\bg(\bm) = \frac{1}{K} \sum_{j=1}^K \nabla_{\bm} (\| \bf{B}\bf{w}_j\|_2^2),
\]
so in particular $\mathbb{E}(\bg)$ is the true gradient of $\phi(\bm)$, while $\bG$ is any desired Hessian approximation 
that depends on the same set of shots as $\bg$. To implement each iteration, we again work with $K$ PDEs; 
the variables $\bf{w}_j$ can vary between iterations. 

These randomized accelerations have a significant impact in full acquisition geometries, where 
all sources are recorded by the all receivers.  
In practical scenarios, this condition nearly always fails  because of budgetary and physical constraints. 
For example, in marine seismic acquisition, both sources and receivers have to move with the ship, 
so all-see-all is impossible by definition. 

We can model any acquisition geometry as a {\it subsample} of a hypothetical all-see-all acquisition.
Call the fully sampled data $\bDf$ and let  
$\bM$ denote the mask that extracts available observations: 
\begin{equation}
\label{eq:mask}
  \bM^{ij}=\left\{
  \begin{array}{@{}ll@{}}
    1, & \mbox{if}\ (i,j) \in \Omega \\
    0, & \mbox{otherwise},
  \end{array}\right.
\end{equation} 
where $\Omega$ is a observed subset of entries from the fully sampled data matrix $\bDf$. 
We then have $\bDs = \bM \odot \bDf$, where $\odot$ is the element-wise (Hadamard) product. 

The mask $\bM$ precludes a straightforward application of the simultaneous shot method. 
The masked version of~(\ref{eq:reduced}) is given by 
\begin{equation}
\label{eq:masked}
\min_{\bm}\frac{1}{2} \mathbb{E}_{\bw}\| \big({\bM\bf{ \odot (P^T H^{-1}({m}) Q )\big) {w_j} - \bDs{w_j}\| ^2_F}}.
\end{equation}
We no longer have the simultaneous source term, 
because  we cannot  simply use the matrix vector product  ${\bf{Q w_j}}$. 
Instead, we have to compute ${\bM\bf{\odot (P^T H^{-1}({m}) Q }})$ first before we can multiply it with ${\bf{w_j}}$. 
This negates the entire motivation of the simultaneous shot method, since we have to solve all PDEs at every 
iteration. 

Our approach is to complete the data, so that we can work with $\bDf$ even though we never observed it 
in the first place. We discuss data-driven low-rank interpolation techniques in the next section, 
and then present a unified interpolation and inversion approach.

\section{Using Low-Rank Interpolation as a Pre-processing Step}
\label{sec:lowrank}

In this section, we review optimization-based methods for interpolating missing traces 
to reconstruct fully sampled data $\bDf$ from the partial observations $\bDs$. 
Commonly used interpolation techniques promote a parsimonious representation 
of the data in a transform domain. For example, if the vectorized data is {\it compressible} 
in a particular domain, we interpolate by penalizing a sparsifying penalty
(such as the 1-norm) on the coefficients of $\bDf$ in that domain; Fourier~\cite{sacchi1998interpolation}, Wavelet, 
and Curvelet~\cite{herrmann2008non} domains are frequently used. 
Analogously, if the full data can be organized into a 2-dimensional matrix with quickly decaying singular values, we can look for low-rank decompositions that match available observations~\cite{oropeza2011simultaneous,kumar2015efficient,aravkin2014fast}.
Here, we focus on low-rank approaches, which are computationally and memory efficient for large-scale seismic data problems~\cite{kumar2015efficient}.

To recover fully sampled data from the subsampled data, we can solve the following rank-minimization problem

\begin{eqnarray}
& \min_{\bDf} \qquad \qquad\mbox{rank}(\bDf) \nonumber\\
& \mbox{subject to} \qquad \frac{1}{2}\| \bM \odot \bDf - \bDs\|_F^2 \leq \epsilon\label{eq:rankmin},
\end{eqnarray}
where $\epsilon$ specifies how closely entries of $\bM\odot \bDf$ must be to the actual observed 
entries $\bDs$. Problem~(\ref{eq:rankmin}) is NP hard, and algorithms that provide exact solutions
have complexity that is doubly exponential in the dimension of the matrix~\cite{candes2009exact,candes2010matrix}. 
Popular alternatives are (1) a convex relaxation using  the nuclear-norm, 
which replaces $\mbox{rank}(\bDf)$ by the sum of singular values of 
$\bDf$~\cite{candes2009exact,candes2010matrix}, or (2) 
an explicit factorization $\bDf = \bL \bR$, where the modeler selects an upper bound $k$
for the rank of the factors a priori~\cite{kumar2015efficient,aravkin2014fast} and solves 
\begin{eqnarray}
& \min_{\bL, \bR}  \qquad \frac{1}{2}\|\bL\|_F^2 + \frac{1}{2} \|\bR\|_F^2\\
& \mbox{subject to} \qquad \frac{1}{2}\| \bM \odot (\bL\bR) - \bDs\|_F^2 \leq \epsilon\label{eq:rankfact}.
\end{eqnarray}
The regularizer is an upper bound to the nuclear norm (sum of singular values) of $\|\bL \bR\|_*$~(\cite{aravkin2014fast}):
\[
\|\bL\bR\|_* \leq  \frac{1}{2}\|\bL\|_F^2 + \frac{1}{2} \|\bR\|_F^2.
\]

Once we recover $\bDf$, we can then apply stochastic optimization or solve an SAA approximation, 
as described in the introduction, without concern for the acquisition-encoding mask $\bM$.

\subsection{Numerical Example} \label{Coupled}

\begin{figure}[h!]
\centering
\includegraphics[width=5in]{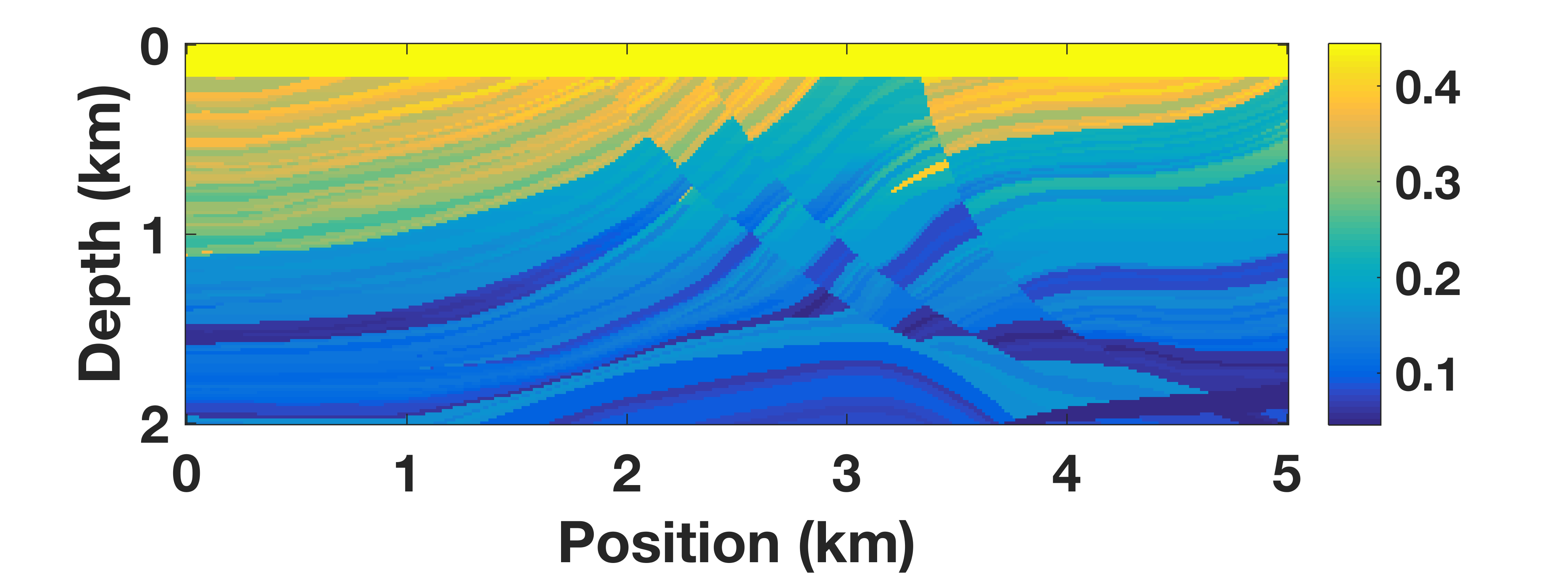}
\includegraphics[width=5in]{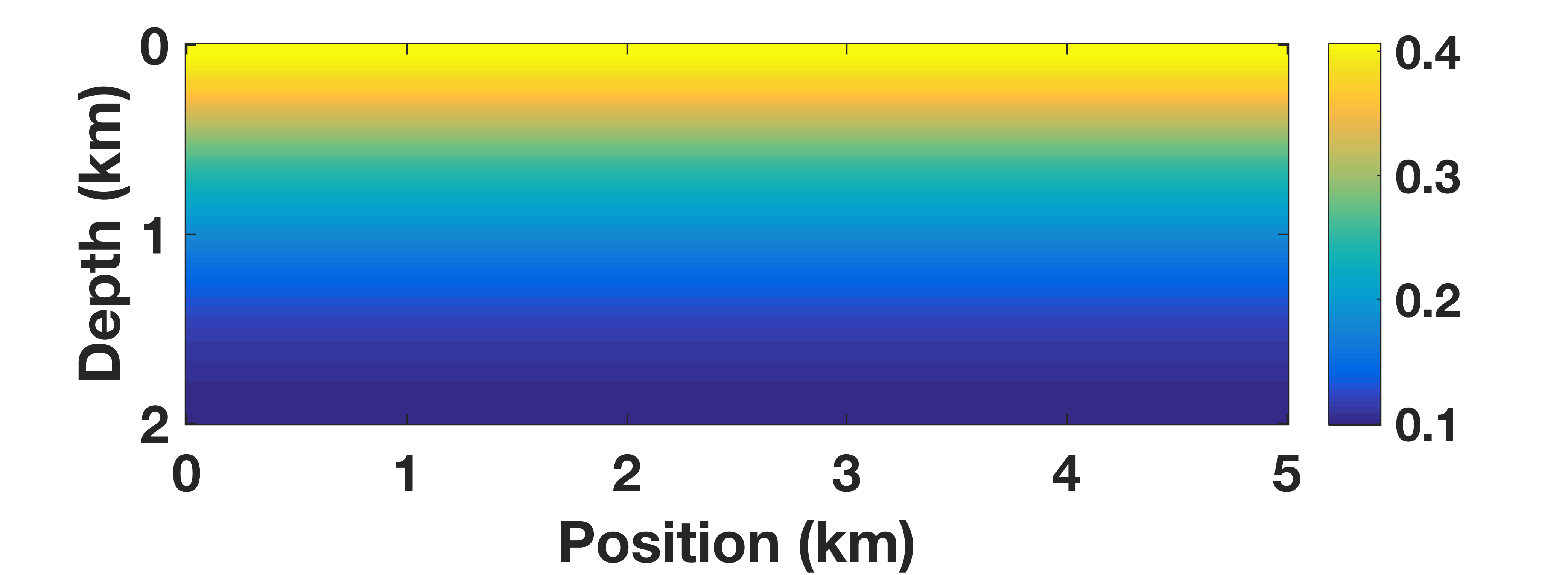}
\caption{\label{Truemod}
\emph{Left}: True Marmousi velocity model. \emph{Right}: Initial velocity model used in waveform-inversion .}
\end{figure}

In this section, we apply stage-wise analysis (interpolation followed by inversion) to 
solve a full-waveform inverse problem over a complex acquisition geometry. 
Here, $\bH$ represents the constant-density acoustic Helmholtz wave equation
\[
\bH =  \omega^2 {\bf{m}}^2 + {\nabla}^2,
\] 
and $\omega$ encodes temporal frequencies and $\bf{m}$ represents the squared-slowness. 
We simulate data using the full Marmousi velocity model \cite{brougois1990marmousi}, 
which has complex geological structure with steeply dipping events (Figure~\ref{Truemod}).
We use a model grid spacing of 15 m and simulate a fixed-spread acquisition configuration, 400 co-located sources and receivers with 12.5 m spacing, 
and a Ricker source wavelet with a peak frequency of 15 Hz. 
To generate a dataset with partial observations, 
we use a frequency spectrum of 3-30 Hz and interpolate each frequency slice independently in this range using rank-minimization framework.

\begin{figure}[h!]
\centering
\includegraphics[width=5in]{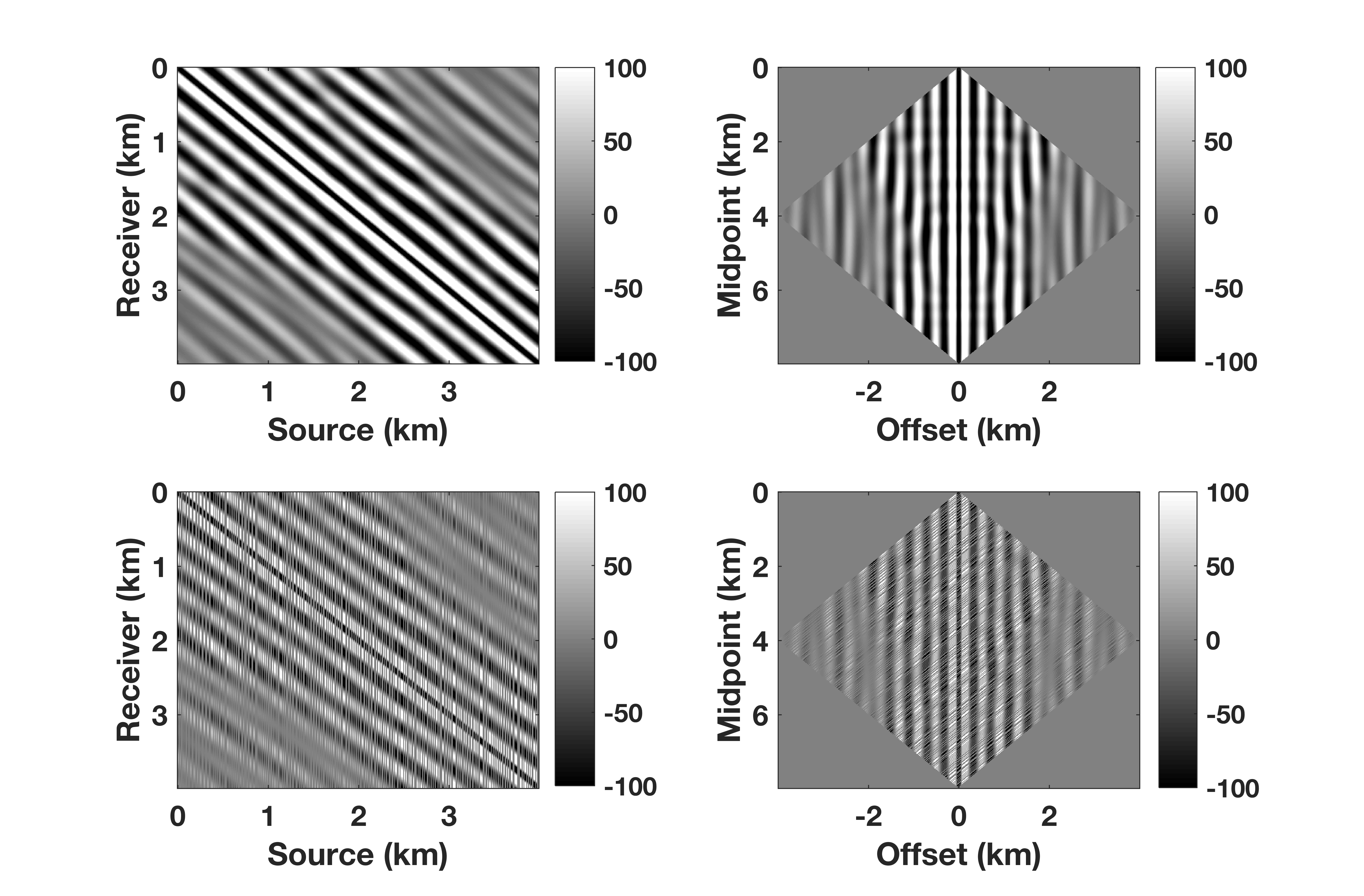}
\caption{A frequency slice at 4 HZ from the Marmousi model. \emph{Top}: Fully sampled data. \emph{Bottom} 50\% subsampled data. \emph{Left}: Source-Receiver domain. \emph{Right}: Midpoint-offset domain.}
\label{MHplot}
\end{figure}

\begin{figure}[h!]
\centering
\includegraphics[width=3in]{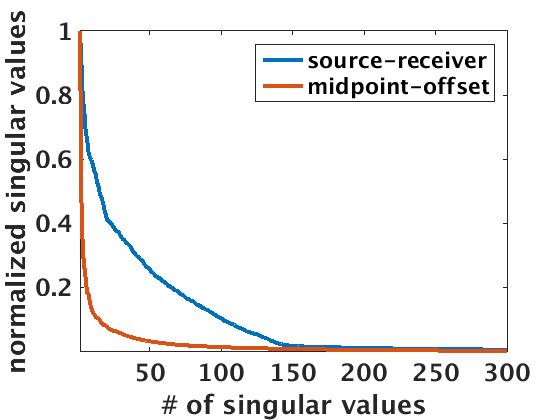}
\includegraphics[width=3in]{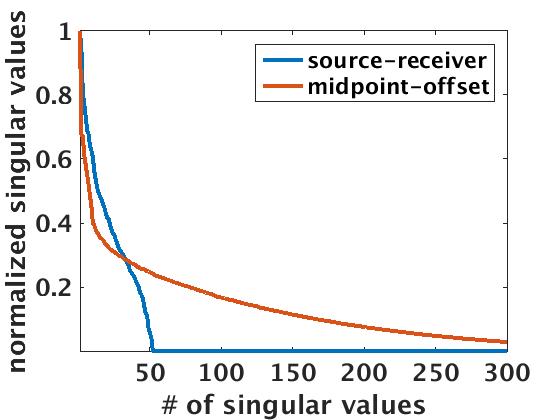}
\caption{\label{SVDplot4Hz}
SVD plot of the frequency slice at 4 HZ from the Marmousi model. \emph{Left}: SVD of fully sampled data. \emph{Right} SVD of 80\% subsampled data.}
\end{figure}

In order to recover missing entries using low-rank optimization, we require that: 
\begin{enumerate}
\item The fully sampled data should have quickly decaying singular values. 
\item Subsampling the data should destroy this fast decay. 
\end{enumerate}
When these conditions are met, we can penalize a rank proxy in order to recover the full data volume.
For seismic data, monochromatic frequency slices satisfy these requirements~\cite{aravkin2014fast} 
in the midpoint-offset domain, where the midpoint ($m_i$) and offset ($h_i$) are defined as
 \begin{equation}
\label{eq:midoff}
m_i = \frac{r_i + c_i}{2}, \quad h_i = \frac{r_i - c_i}{2}, \quad \mathcal{T}: (r_i, c_j) \rightarrow (m_i, h_i),
\end{equation}
with $r_i$ and $c_j$ the  $i^{th}$ row (source position) and $j^{th}$ column (receiver position) respectively.
The transformation $\mathcal{T}$ rotates the data matrix by 45$^{\circ}$ clockwise as illustrated in Figure \ref{MHplot}.  
Figure~\ref{SVDplot4Hz} compares the singular value decay of a monochromatic seismic data slice at 4 Hz in the source-receiver and midpoint-offset domain; conditions (1) and (2) above are satisfied in the midpoint-offset domain, but not in the source-receiver domain (Figure \ref{SVDplot4Hz}). Therefore, we formulate a rank-minimization problem~(\ref{eq:rankfact}) in the midpoint-offset domain:

\begin{eqnarray}
& \min_{\bL, \bR}  \qquad \frac{1}{2}\|\bL\|_F^2 + \frac{1}{2} \|\bR\|_F^2\nonumber\\
& \mbox{subject to} \qquad \frac{1}{2}\| \bM \odot (\mathcal{T^*}(\bL\bR)) - \bDs\|_F^2 \leq \epsilon\label{eq:rankmhfact}.
\end{eqnarray}

\begin{figure}%
\centering
\subfigure[][]{%
\label{fig:ex3-b}%
\includegraphics[width=5in]{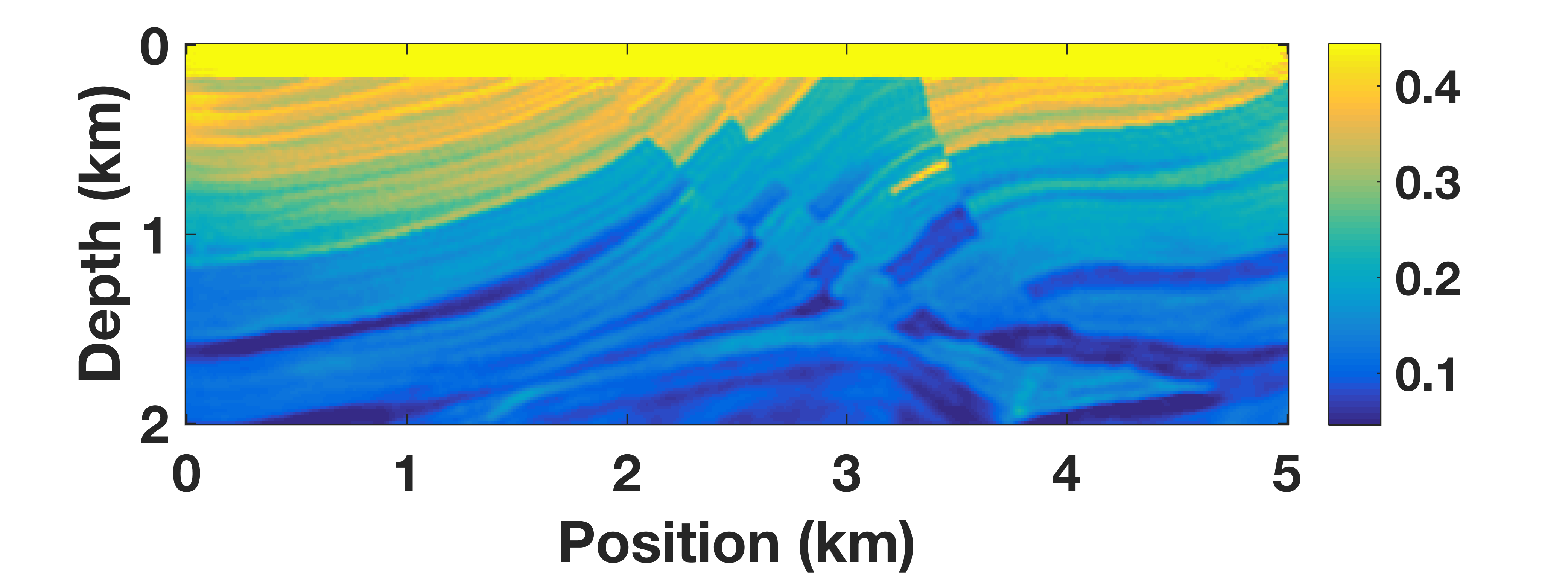}} \\
\subfigure[][]{%
\label{fig:ex3-c}%
\includegraphics[width=5in]{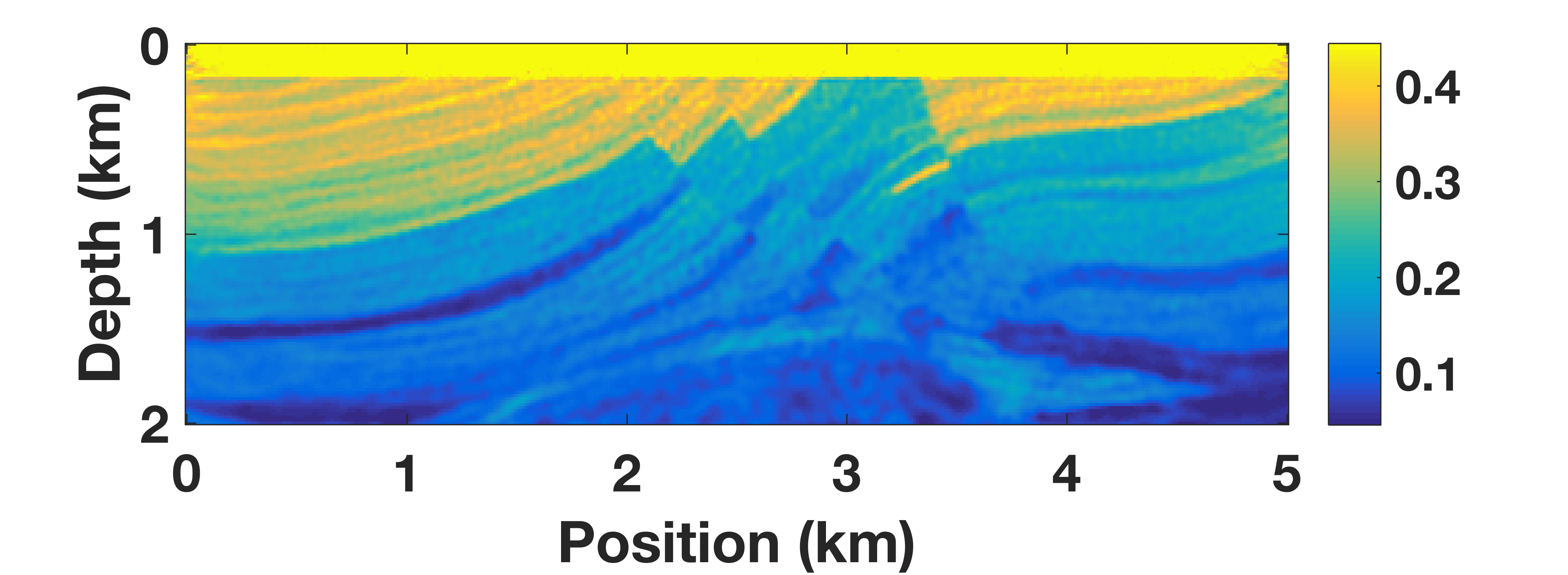}} \\
\subfigure[][]{%
\label{fig:ex3-d}%
\includegraphics[width=5in]{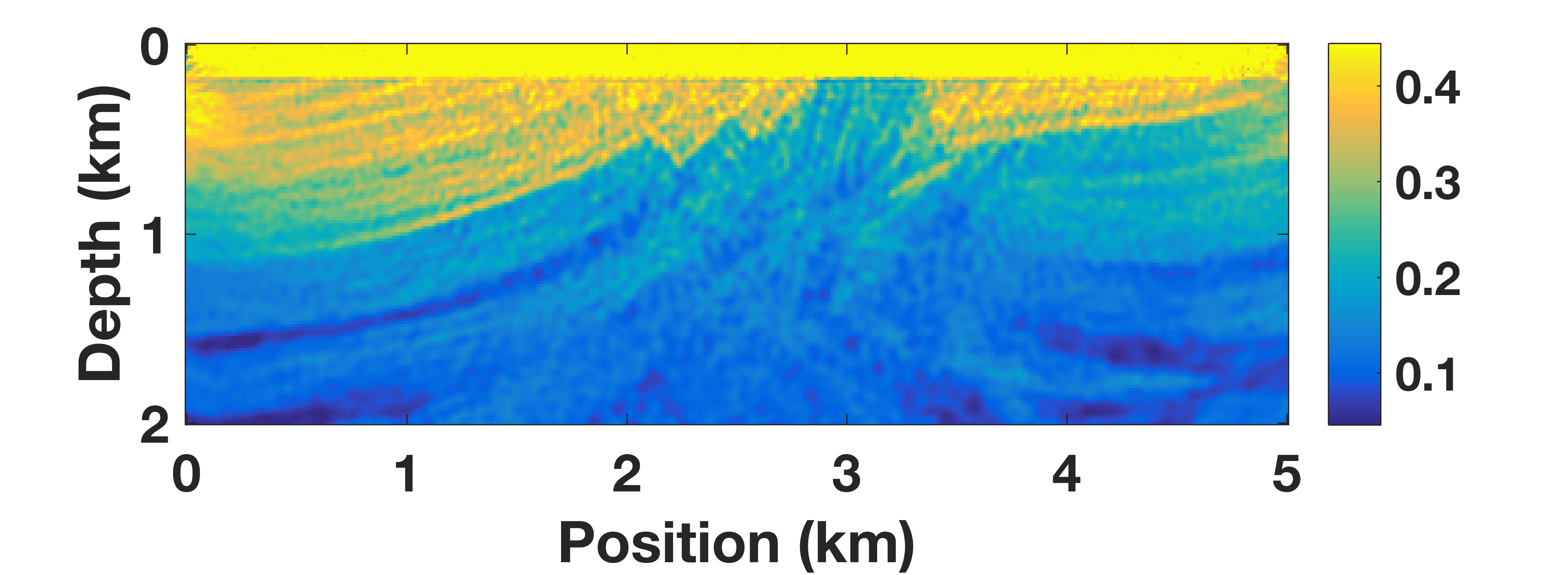}}%
\caption[A set of three subfigures]{\label{disjointModels}
A set of three subfigures  computed using the disjoint inversion approach with 25 PDE solves:
%\subref{fig:ex3-a} the initial model;
\subref{fig:ex3-b} recovered model from 50\% missing data;
\subref{fig:ex3-c} recovered model from 75\% missing data; and,
\subref{fig:ex3-d} recovered model from 85\% missing data.}%
\label{fig:ex3}%
\end{figure}
 
  \begin{figure}[h!]
\centering
\includegraphics[width=5in]{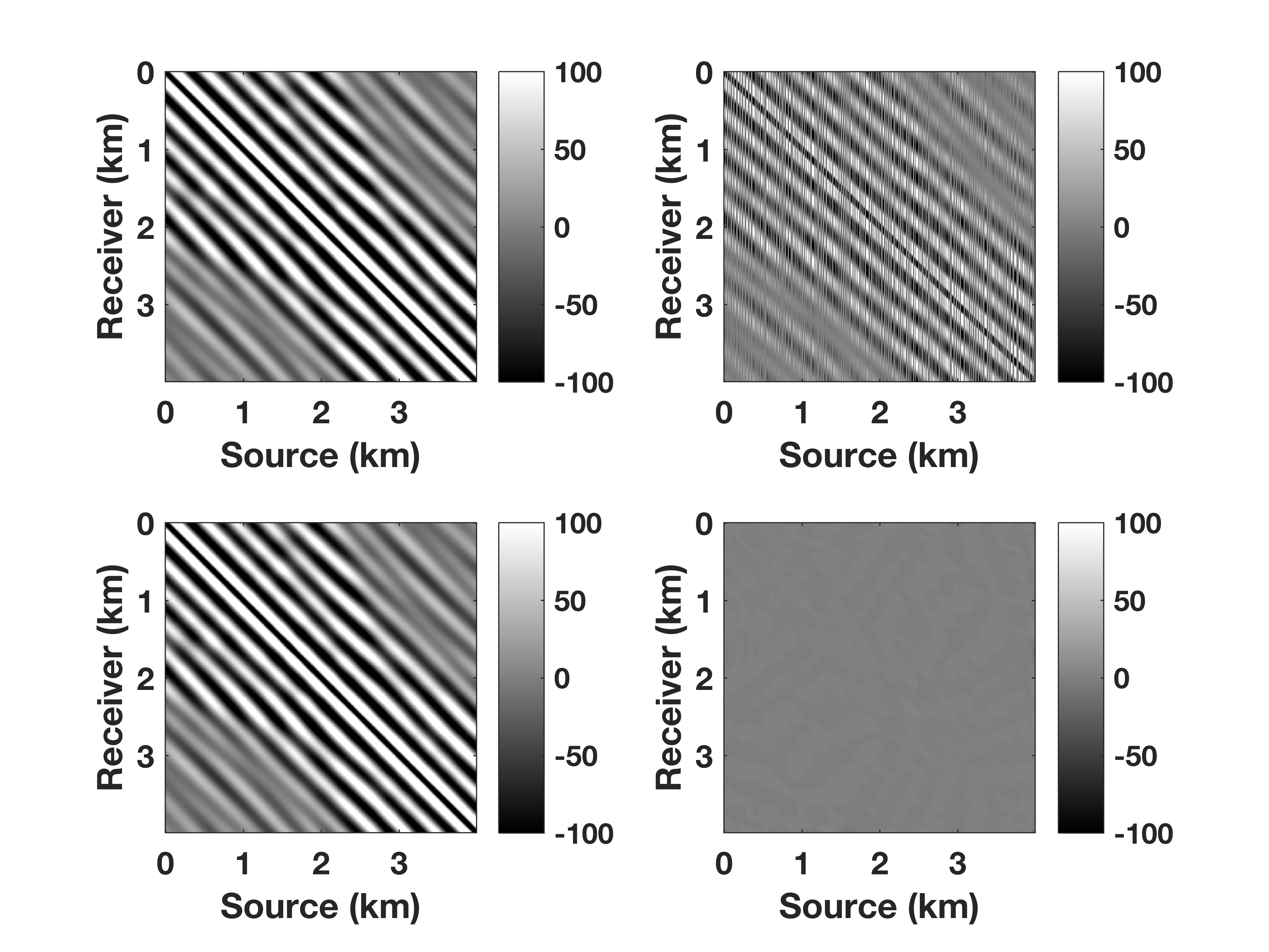}
\caption{\label{Data50_4Hz}
A frequency slice at 4 HZ from the Marmousi model in the acquisition (i.e. source-receiver) domain. Signal to noise ratio is 27.3105 dB. \emph{Top Left}: True Data. \emph{Top Right} 50\% subsampled data. \emph{Bottom Left}: Recovered data. \emph{Bottom Right}: Difference between true and recovered data }
\end{figure}

 \begin{figure}[h!]
\centering
\includegraphics[width=5in]{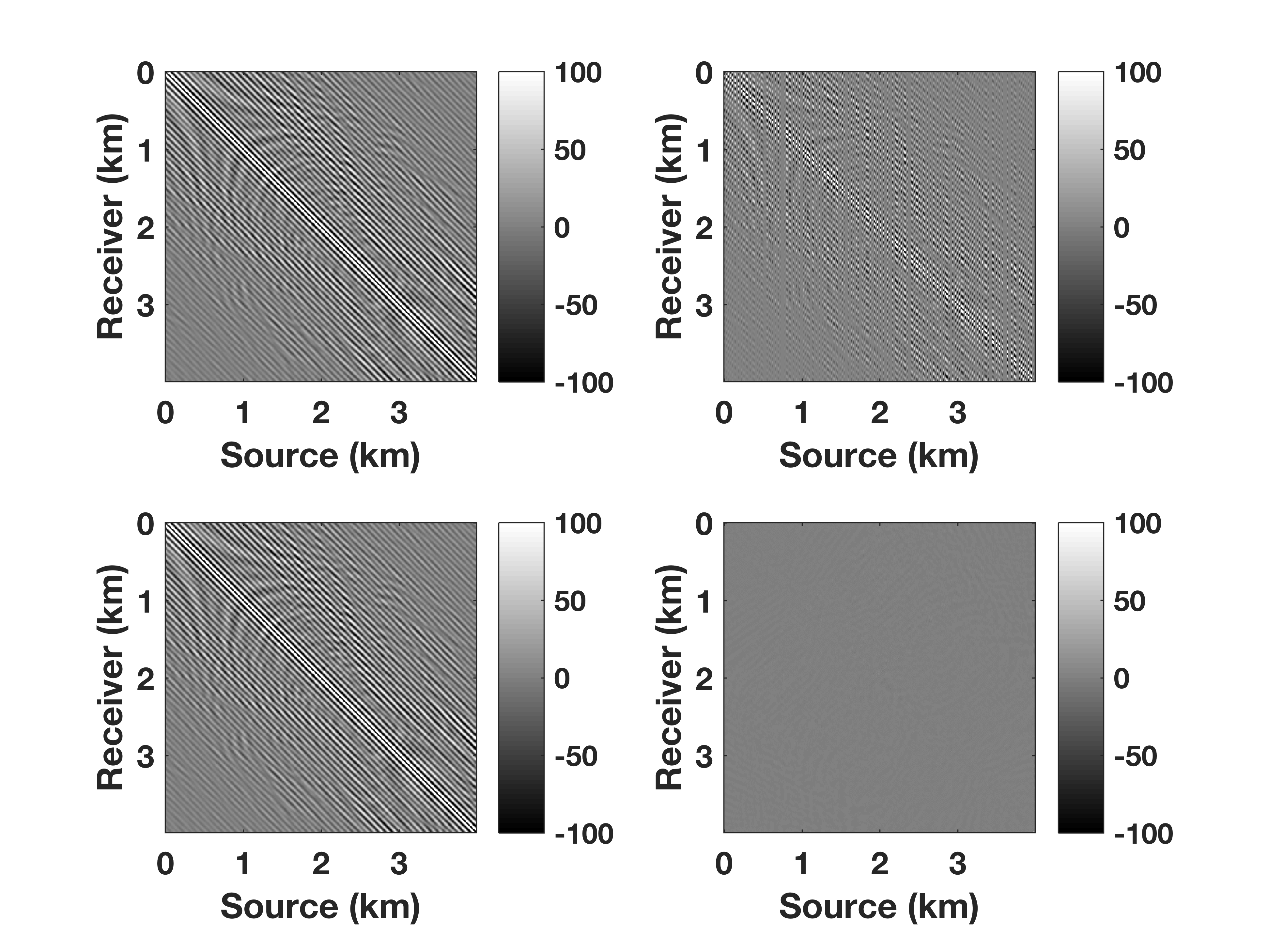}
\caption{\label{Data50_18Hz}
A frequency slice at 18 HZ from the Marmousi model in the acquisition (i.e. source-receiver) domain. Signal to noise ratio is 22.3444 dB. \emph{Top Left}: True Data. \emph{Top Right} 50\% subsampled data. \emph{Bottom Left}: Recovered data. \emph{Bottom Right}: Difference between true and recovered data }
\end{figure}

\begin{figure}[h!]
\centering
\includegraphics[width=5in]{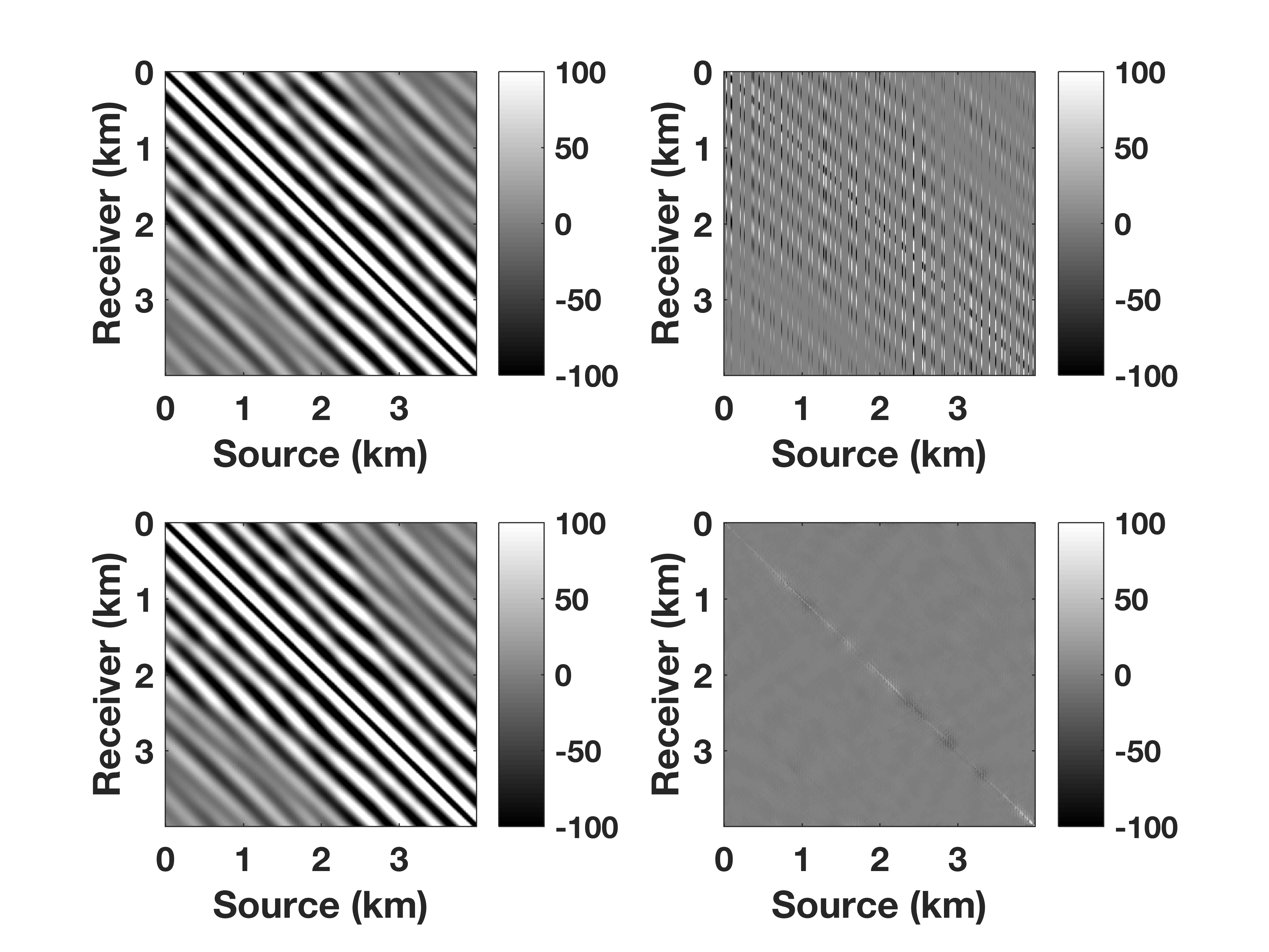}
\caption{\label{Data_4Hz}
A frequency slice at 4 HZ from the Marmousi model in the acquisition (i.e. source-receiver) domain. Signal to noise ration is 21.42 dB. \emph{Top Left}: True Data. \emph{Top Right} 85\% subsampled data. \emph{Bottom Left}: Recovered data. \emph{Bottom Right}: Difference between true and recovered data }
\end{figure}

\begin{figure}[h!]
\centering
\includegraphics[width=5in]{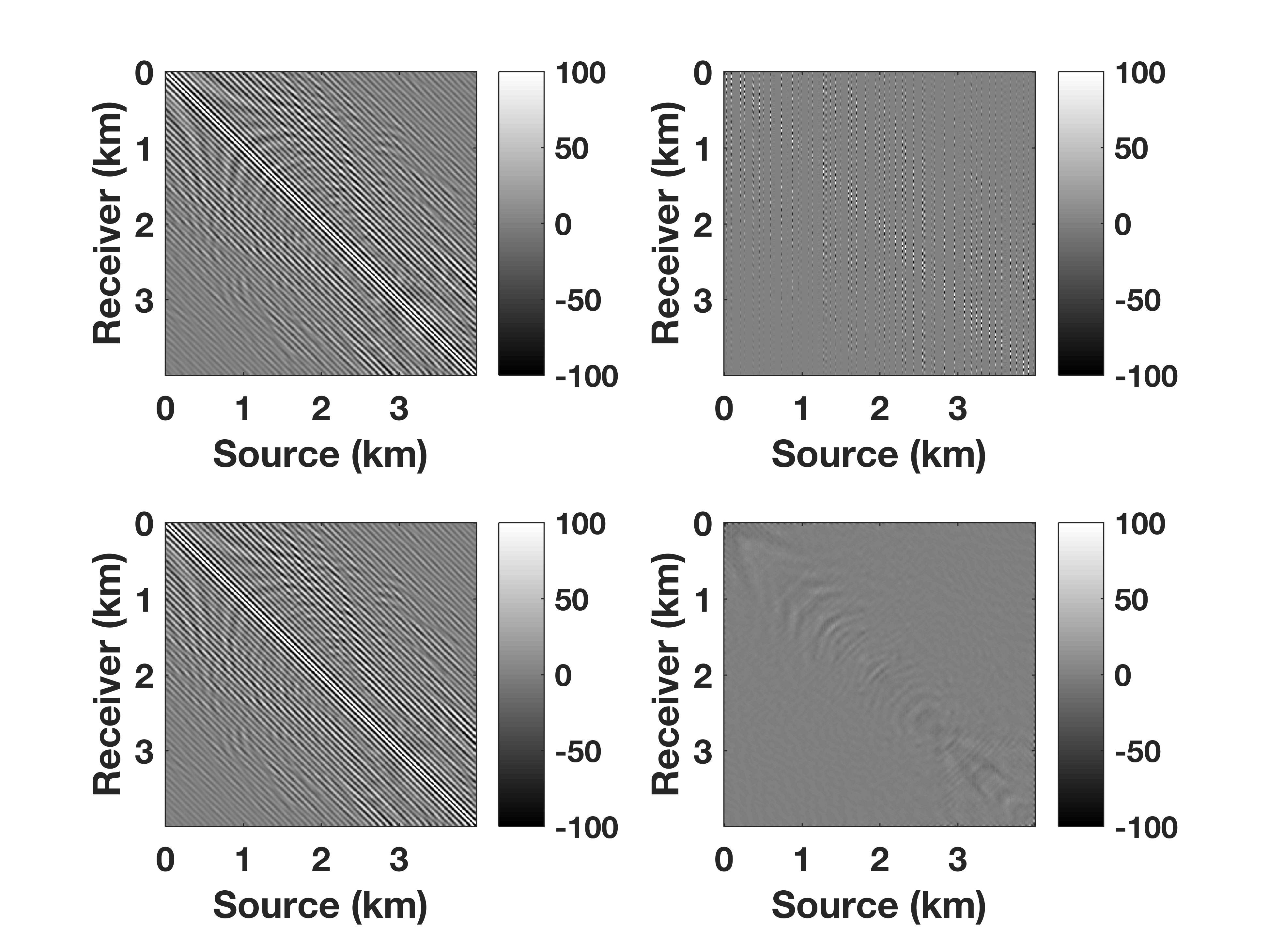}
\caption{\label{Data85_18Hz}A frequency slice at 18 HZ from the Marmousi mode in the acquisition (i.e. source-receiver) domain. Signal to noise ratio is 8.7117 dB. \emph{Top Left}: True Data. \emph{Top Right} 85\% subsampled data. \emph{Bottom Left}: Recovered data. \emph{Bottom Right}: Difference between true and recovered data }
\end{figure}

Once we solve~(\ref{eq:rankmhfact}) for $\bDf$ in the midpoint-offset domain, we 
apply  $\mathcal{T^*}$ to the recovered data to map it back to the source-receiver domain. 
We then solve waveform inversion over an all-see-all geometry using simultaneous shots. 
 We use a linear gradient model to initialize, as shown in Figure \ref{Truemod}. 
 The inversions are carried out sequentially in ten overlapping frequency bands on the interval 3---30 Hz (\cite{bunks1995multiscale}), each using 25 different randomly selected simultaneous shots and six selected frequencies in each band with an interval of 0.2 Hz. 
 We use a limited-memory L-BFGS method~\cite{wright1999numerical}. 
 
Figure \ref{disjointModels} shows the inversion results with interpolated data for 50\%, 75\% and 85\% missing entries. We get excellent inversion results for 50\% missing data, since we get good reconstruction quality data using rank-minimization based framework (see Figures \ref{Data50_4Hz} and \ref{Data50_18Hz}). We start to see deterioration of  inversion results at about 75\% missing data.  Finally, as we move to 85\% missing data, inversion quality deteriorates noticeably. This is due to the fact that the data reconstruction using rank-minimization is poor as we move from low to high subsampling ratio (see Figures \ref{Data_4Hz} and \ref{Data85_18Hz}).

We can still recover a reasonable solution by solving the waveform inversion problem using the 
15\% of the observed data by applying classical approaches (non-randomized) FWI techniques. 
However, this is a costly prospect, as we must solve all PDEs in each iteration to get the predicted data at the available locations.  Simultaneous shots can get the solution quickly by solving very few PDEs using simultaneous sources, but require fill-in for unseen data, which introduces error and degrades the solution quality.
There appears to be a clear speed vs. quality tradeoff.  In the next section we show that we can get around 
this bottleneck by using a unified approach, where we combine low-rank interpolation with inversion. 

\section{Simultaneous Data Completion and Inversion} \label{sec:unified}
From the previous section, it is clear that while filling in an all-see-all acquisition scenario 
makes it possible to use simultaneous shots, it also introduces error when the available data 
comprise less than $25\%$ of the hypothetical full volume. 
To push past this boundary, we propose a unified data completion and waveform inversion approach, 
and develop a customized algorithm to solve it. 

The approach is intuitive: we merge the low-rank data fill-in problem~(\ref{eq:rankmhfact})
with the waveform inversion problem~(\ref{eq:reduced}). The joint problem can be formulated as follows:
\begin{eqnarray}
\min_{\bL, \bR, \bm} &  \qquad \frac{1}{2} \|\bL\|_F^2 + \frac{1}{2} \|\bR\|_F^2\label{eq:jointOpt1}\\
\mbox{subject to} & \qquad  \| \bM \odot (\mathcal{T}^*(\bL\bR)) - \bDs \|_F^2 + 
\frac{\lambda}{2}\| \bP \bH(\bm)^{-1} \bQ- (\mathcal{T}^*(\bL\bR)) \|_F^2 \leq \epsilon\nonumber, 
\end{eqnarray}
where $\lambda$ is a tradeoff parameter. The advantage of~(\ref{eq:jointOpt1}) over the stage-wise approach is that the 
model $\bm$ informs the interpolation. The data-fit term in the constraints 
prevents fill-in that is inconsistent with the physics.  The challenge of~(\ref{eq:jointOpt1}) is the need to solve $\bP \bH(\bm)^{-1} \bQ$ for all the sources, which is prohibitively expensive for large-scale seismic data. To mitigate this, we design a stochastic block-coordinate descent algorithm for~(\ref{eq:jointOpt1}), where we need to solve only a few PDEs in each subproblem by using simultaneous shots. We can encode simultaneous shots explicitly using a random matrix 
$\bW\in \mathbb{R}^{(N_s \times K)}$, where each column encodes a simultaneous shot by drawing a random vector from a standard i.i.d Gaussian distribution. The shots $\bW$ are fixed throughout each iteration of the block coordinate 
descent method described below; after each iteration described below is complete, $\bW$ is resampled.

\noindent
{\bf Solving for $\bL, \bR$ for fixed $\bm$}.
For fixed $\bm$, we consider a problem~(\ref{eq:jointOpt1})
is a residual-constrained matrix factorization problem, studied by~\cite{aravkin2014fast}. 
The randomized subproblem using the simultaneous shots $\bW$ is given by 
\begin{eqnarray}
\min_{\bL, \bR} &  \qquad \frac{1}{2} \|\bL\|_F^2 + \frac{1}{2} \|\bR\|_F^2\label{eq:rBPDN}\\
\mbox{subject to} & \qquad  \| \bM \odot (\mathcal{T}^*(\bL\bR)) - \bDs \|_F^2 + 
\frac{\lambda}{2}\|  \bF \bW- (\mathcal{T}^*(\bL\bR)) \bW \|_F^2 \leq \epsilon\nonumber,
\end{eqnarray}
where $\bF \bW =  \bP \bH(\bm)^{-1} \bQ \bW$ is the forward-modeled simultaneous seismic data at the current $\bm$.

We solve this problem using the extension of the SPG$\ell_1$ solver~\cite{van2008probing}
developed by~\cite{aravkin2014fast}\footnote{https://github.com/UW-AMO/spgl1} 
to handle matrix factorization. The core idea is to solve a sequence of Lasso subproblems:
\begin{eqnarray}
v(\tau) := \min_{\bL, \bR} &  \qquad \| \bM \odot (\mathcal{T}^*(\bL\bR)) - \bDs \|_F^2 + 
\frac{\lambda}{2}\|  \bF \bW- (\mathcal{T}^*(\bL\bR)) \bW \|_F^2\nonumber\\
\mbox{subject to} & \qquad  \frac{1}{2} \|\bL\|_F^2 + \frac{1}{2} \|\bR\|_F^2 \leq \tau\label{BCU_lasso}, 
\end{eqnarray}
where the $\tau$ corresponding to $\epsilon$ is found by solving $v(\tau) = \epsilon$ with
a root finding method~\cite{aravkin2016level}. The subproblem~(\ref{BCU_lasso}) is solved using a 
fast projected gradient method. The gradient computations for~(\ref{BCU_lasso}) are straightforward. 
The term $\bF \bW$ needs to be computed only at the beginning of the interpolation for the current $\bm$. We  also reuse this precomputed $\bF \bW$ to evaluate the gradient at the first iteration of the $\bm$ subproblem described below.

\noindent
{\bf Solving for $\bm$ with fixed $\bL, \bR$.} When $\bL, \bR$ are fixed, problem~(\ref{eq:jointOpt1})
is essentially a feasibility problem in $\bm$. For a set of simultaneous shots $\bW$, we want to solve  
\begin{equation}
\label{eq:rFWI}
\min_{\bm} \frac{1}{2}\| \bP \bH(\bm)^{-1} \bQ \bW - \mathcal{T}^{*}(\bL\bR) \bW \|_F^2.
\end{equation}
$\bL\bR$ is fixed for this subproblem, and we simply use L-BFGS to update $\bm$. 
L-BFGS only needs gradient information, and the gradients are computed 
using the standard adjoint-state method~\cite{plessix2006review}. 

The overall algorithm is summarized in Algorithm~\ref{alg:unified}. 
The subproblems for $\bL, \bR$ and $\bm$ terminate when an iteration cap is reached. 
\begin{algorithm}[H]
\caption{Unified Algorithm for Fast Interpolated FWI}
\label{alg:unified}
\begin{algorithmic}[1]
\State{\bfseries Input:} $\bM$, $\bDs$, $\mathcal{T}$, $\bP$, $\bQ$
\State{\bfseries Initialize:} $k = 1$, $\bL^1$, $\bR^1$, $\bm^1$

\While{not converged}
\State{Sample} $\bW^k$ so that $\mathbb{E}(\bW^k(\bW^k)^T) = I$
\Let{$\bL^{k+1}, \bR^{k+1}$}{$\arg\min_{\bL, \bR}$(\ref{eq:rBPDN}) with $\bm^k, \bW^k$}
\Let{$\bm^{k+1}$}{$\arg\min_{\bm}$~(\ref{eq:rFWI}) with $\bL^{k+1}, \bR^{k+1}, \bW^k$}
\Let{$k$}{$k+1$}
\EndWhile
\State{\bfseries Output:} $\bL^k, \bR^k, \bw^k$
\end{algorithmic}
\end{algorithm}

\section{Numerical Experiment}
We test the joint inversion approach on the same dataset as the disjoint inversion, 
and compare the results with those obtained previously by the stage-wise approach. 
We initialize L and R following~\cite{aravkin2014fast} by using the SVD of 
the subsampled data in the midpoint offset domain:
\[
[{\bf U,S,V}] = SVD(\bDs),
\] 
and set $\bL = {\bf U\sqrt{S}}$ and $\bR = {\bf V\sqrt{S}}$ as the initial value. 
The models were recovered by performing 25 partial solves 
of the $\bm$ and $\bL, \bR$ subproblems. 
%iterations in the main loop of~\cref{BCUM}. 
Instead of using the full 400 PDE solves to recover the model, 
we needed only 25 PDEs, just as for the disjoint inversion experiment.

\begin{figure}
\centering
\subfigure[][recovered model from 75\% missing data]{%
\label{fig:ex4-c}%
\includegraphics[width=5in]{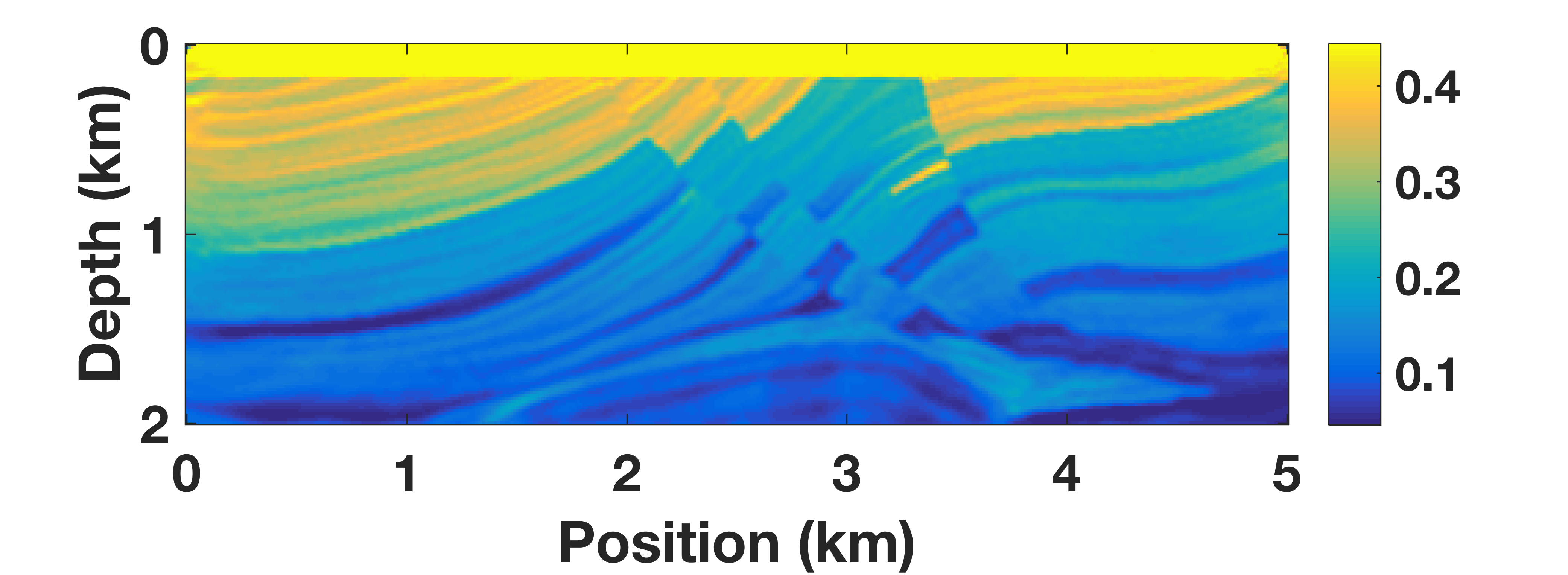}}%
\hspace{8pt}%
\subfigure[][recovered model from 85\% missing data]{%
\label{fig:ex4-d}%
\includegraphics[width=5in]{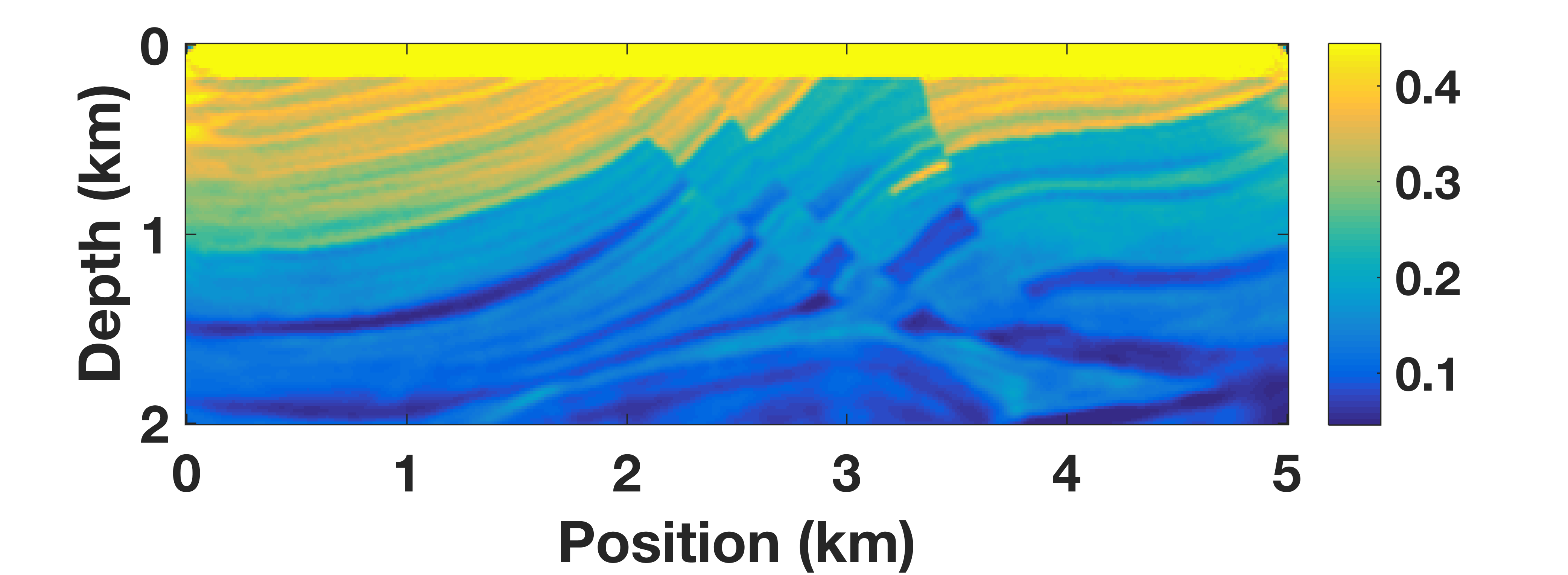}}%
\caption[A set of four subfigures]{\label{JointModels}Results computed using the joint inversion approach with 25 PDE solves.
%\subref{fig:ex4-a},
%\subref{fig:ex4-b},
%\subref{fig:ex4-c} Results for 75\% missing data, \subref{fig:ex4-d} 85\% missing data.
} 
\label{fig:ex4}%
\end{figure}

\begin{figure}[h!]
\centering
\includegraphics[width=1.0\textwidth]{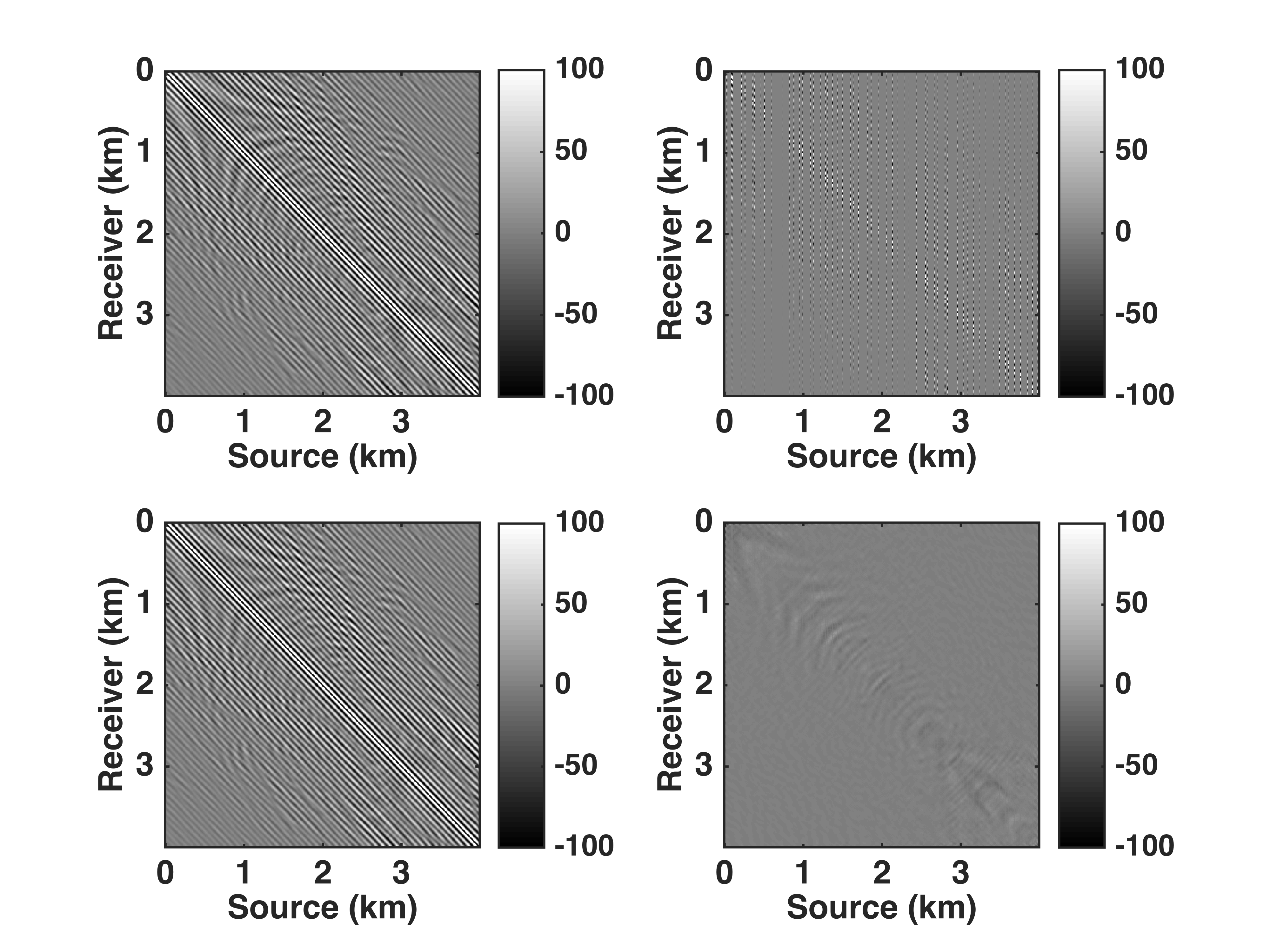}
\caption{\label{Data85JOINT_18Hz}A frequency slice at 18 HZ from the Marmousi model in the acquisition (i.e. source-receiver) domain. Signal to noise ratio is 13.96 dB. \emph{Top Left}: True Data. \emph{Top Right} 85\% subsampled data. \emph{Bottom Left}: Recovered data. \emph{Bottom Right}: Difference between true and recovered data }
\end{figure}

Figure \ref{JointModels} shows the joint inversion results for 75\% and 85\% missing data. 
Stage-wise and unified inversion give the same quality results for 50\% missing data, hence, we do not plot the similar figure for joint inversion. 
However, the results for unified inversion at higher levels of missing data are far better. 
 The recovery of the model is excellent even at 85\% missing data, 
 with good data reconstruction at high frequencies, see Figure \ref{Data85JOINT_18Hz}.

\section{Discussion and Conclusion}
We have presented a joint interpolation and inversion framework, which can recover the high-fidelity seismic data and subsurface velocity model under high-subsampling scenarios (e.g. 80\% missing data). The proposed approach exploits the fast singular value decay of seismic data in the midpoint-offset domain, and allows fast stochastic methods to bear on the inverse problem despite the missing data.  In particular, we recover missing data and then use simultaneous shots from the recovered maps to inform waveform inversion. Using carefully selected stylized examples, we showed that our method outperforms the traditional seismic data reconstruction and inversion methods, is fast, and highly scalable. To our knowledge, this work is first of a kind in combining seismic data interpolation and waveform-inversion in a joint framework, allowing us to efficiently and simultaneously reconstruct seismic data volumes and invert for artifact-free velocity models of the subsurface. 

This paper opens several new avenues of research.  First, so far we assume that the starting velocity model is not cycle skipped, which occurs when predicted and observed data differ by more than half a cycle. For future work, we plan to relax this assumption by using Wavefield Reconstruction Inversion~\cite{vanleeuwen2015IPpmp} to handle the cycle skipping phenomenon. 
Second, here we consider only randomly subsampled data scenarios rather than structured scenarios such as 
towed-streamer marine seismic acquisition. Future work  includes developing the proposed approach for 
these complex geophysical acquisition scenarios.  
Third, Seismic data are typically irregularly sampled along spatial axes,  and we also plan to adapt the approach to non-uniform sampling grids. Finally, we expect that these methods can 
be extended to perform joint interpolation and inversion for large-scale 3-D seismic data acquisition, where the underlying model is 3D and the observed seismic data is 5D.

\section*{Acknowledgments}
We would like to acknowledge the assistance of volunteers in putting together this example manuscript and supplement. R.K.
would like to thank the member organizations of the SINBAD II project and the SINBAD Consortium for supporting this work.
The work of Prof. Aravkin was supported by the Washington Research Foundation Data Science Professorship.

\bibliographystyle{siamplain}
\bibliography{references}

\end{document}